\documentclass[11pt]{article}
\usepackage{epic,latexsym,amssymb}
\usepackage{tikz}

\usepackage{tikz,epsf}

\usepackage{amsfonts}
\usepackage{amscd}
\usepackage{amsmath}
\usepackage{graphicx}
\usepackage{color}
\usepackage{caption,subcaption}

\newtheorem{thm}{Theorem}
\newtheorem{conj}{Conjecture}

\newtheorem{claim}{Claim}

\newcommand{\cL}{\mathcal{L}}
\newcommand{\R}{\mathcal{R}}

\newcommand{\smallqed}{{\tiny ($\Box$)}}
\newcommand{\qed}{$\Box$}

\newcommand{\cart}{\, \Box \,}
\newenvironment{proof}{\noindent\textbf{Proof.}}{\hfill $\square$}

\newcommand{\vertex}{\node[vertex]}
\tikzstyle{vertex}=[circle, draw, inner sep=0pt, minimum size=6pt]

\newcommand{\QEDmark}{\mbox{\textsc{qed}}}
\newcommand{\proofStarter}[1]{\textsc{#1} }

\def\vertex(#1){\put(#1){\circle*{2}}}
\def\vertexo(#1){\put(#1){\circle{2}}}
\def\vert(#1){\put(#1){\circle*{1.5}}}
\def\verto(#1){\put(#1){\circle{1.5}}}
\def\lab(#1)#2{\put(#1){\makebox(0,0)[c]{#2}}}
\definecolor{DarkGreen}{rgb}{0.2, 0.6, 0.3}

\definecolor{electricindigo}{rgb}{0.44, 0.0, 1.0}

\begin{document}

\title{A Vizing-type result for semi-total domination}

\author{
John Asplund\\
Department of Technology and Mathematics\\
Dalton State College \\
Dalton, GA 30720, USA \\
jasplund@daltonstate.edu\\
\\
Randy Davila\\
Department of Pure and Applied Mathematics\\
University of Johannesburg \\
Auckland Park 2006, South Africa \\
%\small {\tt Email: mahenning@uj.ac.za} \\
\\
Department of Mathematics and Statistics \\
University of Houston--Downtown \\
Houston, TX 77002, USA \\
\small {\tt Email: davilar@uhd.edu}\\
\\
Elliot Krop\\
Department of Mathematics \\
Clayton State University \\
Morrow, GA 30260, USA \\
\small{\tt Email: elliotkrop@clayton.edu}
}

\date{}
\maketitle

\begin{abstract}
A set of vertices $S$ in a simple isolate-free graph $G$ is a semi-total dominating set of $G$ if it is a dominating set of $G$ and every vertex of $S$ is within distance 2 or less with another vertex of $S$. The semi-total domination number of $G$, denoted by $\gamma_{t2}(G)$, is the minimum cardinality of a semi-total dominating set of $G$. In this paper we study semi-total domination of Cartesian products of graphs. Our main result establishes that for any graphs $G$ and $H$, $\gamma_{t2}(G\cart H)\ge \frac{1}{3}\gamma_{t2}(G)\gamma_{t2}(H)$.
\end{abstract}


\section{Introduction}
In this paper we study bounds on a recently introduced domination invariant applied to Cartesian products of graphs. At its core, our work is motivated by the longstanding conjecture of V.G.~Vizing \cite{Vizing} on the domination of product graphs, which states that for any graphs $G$ and $H$, $\gamma(G\cart H)\ge \gamma(G)\gamma(H)$. Here, $\gamma(G)$ is the domination number of $G$, which is the minimum size of a set $D$ of vertices so that every vertex not in $D$ is adjacent to some vertex in $D$, and $\cart$ is the Cartesian product of graphs. The breakthrough ``double-projection" result of Clark and Suen \cite{CS} gave the first Vizing-type bound of $\gamma(G\cart H)\ge \frac{1}{2}\gamma(G)\gamma(H)$. Recently, Bre\v{s}ar \cite{B} improved this bound to $\gamma(G\cart H)\ge \frac{(2\gamma(G)-\rho(G))\gamma(H)}{3}$, where $\rho(G)$ is the two-packing number of $G$. For more on attempts to solve Vizing's conjecture over more than five decades since it was stated, see the survey \cite{BDGHHKR}.

Over the years, due to the unyielding nature of the conjecture, devotees have used offshoots of the domination number to attempt Vizing-type inequalities, in hopes of better understanding the difficulties of the original problem. For example, Bre\v{s}ar, Henning, and Rall \cite{BHR} defined the paired and rainbow domination numbers, and Henning and Rall \cite{HR} conjectured a Vizing-type inequality for total domination. This last conjecture was proved by Ho \cite{H}, who showed that for any graphs $G$ and $H$, $\gamma_t(G\cart H)\ge \frac{1}{2}\gamma_t(G)\gamma_t(H)$. In this result, $\gamma_t(G)$ is the total domination number of $G$, which is  the minimum size of a set $T$ of vertices so that every vertex of $G$ is adjacent to some vertex in $T$. A sharp example was given in \cite{HR} and the characterization of pairs of graphs attaining equality is an active problem, see \cite{BHKM} and \cite{LH}.

Since the difference between a totally dominating set and a dominating set is that every vertex in a totally dominating set must be adjacent to some other vertex in that set, while this rule does not have to hold in a dominating set, we find it instructive to consider Vizing-type inequalities for domination invariants that share properties with both domination and total domination. That is, we want to consider some domination function in between domination and total domination. Such a function, first investigated by Goddard, Henning, and McPillan \cite{GHM}, is the \emph{semi-total domination number of $G$}, $\gamma_{t2}(G)$, which is the minimum size of a set of vertices $S$ in $G$, so that every vertex of $S$ is of distance at most $2$ to some other vertex of $S$, and every vertex not in $S$ is adjacent to a vertex in $S$. Although introduced only a few years ago, this function has seen much recent attention, see \cite{HM,HM2,HM3,HM4,M,ZSX}. 

Although we cannot prove it, we believe that $\gamma_{t2}(G\cart H)\ge \frac{1}{2}\gamma_{t2}(G)\gamma_{t2}(H)$ for any graphs $G$ and $H$. Our result depends on the method of Clark and Suen \cite{CS} and requires more careful analysis of semi-total dominating sets. We show that for any graphs $G$ and $H$, $\gamma_{t2}(G\cart H)\ge \frac{1}{3}\gamma_{t2}(G)\gamma_{t2}(H)$.

\medskip
\noindent\textbf{Definitions and Notation.} 
For notation and graph terminology, we will typically follow \cite{HY}. Throughout this paper, all graphs will be considered undirected, simple, connected, and finite. Specifically, let $G$ be a graph with vertex set $V = V(G)$ and edge set $E = E(G)$. %The order and size of $G$ will be written $n = |V(G)|$ and $m = |E(G)|$, respectively. 
Two vertices $v,w \in V$ are neighbors, or adjacent, if $vw \in E$. The open neighborhood of $v\in V$, is the set of neighbors of $v$, denoted $N_G(v)$, whereas the closed neighborhood is $N_G[v] = N_G(v) \cup \{v\}$. The open neighborhood of $S\subseteq V$ is the set of all neighbors of vertices in $S$, denoted $N_G(S)$, whereas the closed neighborhood of $S$ is $N_G[S] = N_G(S) \cup S$. %The degree of a vertex $v\in V$, is denoted $d_G(v) = |N_G(v)|$. The maximum and minimum degree of $G$, will be denoted $\Delta = \Delta(G)$ and $\delta = \delta(G)$, respectively.
The distance between two vertices $v,w\in V$ is the length of a shortest $(v,w)$-path in $G$, and is denoted by $d_G(v,w)$. The \emph{Cartesian product} of two graphs $G(V_1,E_1)$ and $H(V_2,E_2)$, denoted by $G \cart H$, is a graph with vertex set $V_1 \times V_2$ and edge set $E(G \cart H) = \{((u_1,v_1),(u_2,v_2)) : v_1=v_2 \mbox{ and } (u_1,u_2) \in E_1, \mbox{ or } u_1 = u_2 \mbox{ and } (v_1,v_2) \in E_2\}$.

A subset of vertices $S\subseteq V(G)$ is called a \emph{semi-total dominating set} if $N[S]=V(G)$ and for any vertex $u\in S$, there exists a vertex $v\in S$ so that $d(u,v)\le 2$. The \emph{semi-total domination number of $G$}, written $\gamma_{t2}(G)$, is the size of a minimum semi-total dominating set of $G$. A \emph{$2$-packing} is a subset of vertices $T$ of $G$ so that every pair of vertices in $T$ is of distance at least $3$. The size of a maximum $2$-packing of $G$ is called the \emph{$2$-packing number}, and is written $\rho(G)$.

We will also make use the standard notation $[k] = \{1,\ldots,k\}$, and for two vertices $u,v$, we write $u\sim v$ to indicate that $u$ is adjacent to $v$.

%%%%%%%%%%%%%%%%%%%%%%%%%%%%%%%%%%

\section{Main Results}
In this section we provide our main results. We begin by establishing a Vizing's-type result which makes use of the 2-packing number. 
\begin{thm}
For any isolate-free graphs $G$ and $H$, 
\begin{equation*}
\gamma_{t2}(G\:\Box\: H)\geq \rho(G)\gamma_{t2}(H).
\end{equation*}
\end{thm}

\proof \;Without loss of generality, we assume $G$ is a $(\rho, \gamma)$-graph (where $\rho=\gamma$), and let $\{v_1,...,v_{\rho(G)}\}$ be a maximum 2-packing of $G$. Since each vertex from our packing is distance at least 3 from any other vertex of our packing, we observe that for $i=1,...,\rho(G)$, the closed neighborhoods $N[v_i]$ are pairwise disjoint. Let $\{V_1,...,V_{\rho(G)}\}$ be a partition of $V(G)$ such that $N[v_i]\subseteq V_i$. Let $D$ be a $\gamma_{t2}(G\:\Box\: H)$-set. For $i=1,...,\rho(G)$, let $D_i=D\cap (V_i\times V(H))$, and let $H_i=\{v_i\}\times V(H)$. Further, let $S_i$ be a minimum set of vertices in $G\:\Box\: H$ that semi-totally dominates $H_i$, and contains as many vertices in $H_i$ as possible. Then, $S_i\subseteq V_i\times V(H)$. Next suppose that $S_i$ contains a vertex $x$ such that $x$ is not in $H_i$. Then, $x$ is the unique vertex which semi-totally dominates $x'$, for some $x'\in H_i$. Since $x'$ has neighbors, all of which are dominates by vertices in $S_i$, if we replace $x$ by $x'$ in $S_i$, we see that $S_i$ is still semi-total dominating (Since $x'$ is at distance at least 2 from a vertex which dominates one of its neighbors). Moreover, we have found a set of vertices from $G\:\Box H$ that semi-totally dominates $H_i$ and contains more vertices in $H_i$ than does $S_i$, a contradiction. Hence, we have $S_i\subseteq H_i$, and so $S_i$ is a semi-total dominating set of the copy of $H$ in $G\:\Box\: H$ induced by the set $H_i$. Since $D_i$ semi-totally dominates $\{v_i\}\times V(H)$, $|D_i|\geq |S_i|$. Thus,
\begin{equation*}
\gamma_{t2}(G\:\Box\: H)\geq \sum_{i=1}^{\rho(G)}|S_i|\geq \sum_{i=1}^{\rho(G)}\gamma_{t2}(H)=\rho(G)\gamma_{t2}(H). 
\end{equation*}\qed

Next, we prove a Vizing's type result which relies only on the semi-total domination number.  We do this by partitioning minimum semi-total dominating sets into parts that are and are not totally dominating. Notice that for any graph $G$, if $U=\{u_1,\dots, u_k\}$ is a minimum semi-total dominating set of $G$, then $U$ can be separated into two sets, $X$ and $Y$, where $X$ is the set of vertices of $U$ which are adjacent to at least one other vertex of $U$, and $Y=U\setminus X$. We call such sets $X$, \emph{allied} and such sets $Y$, \emph{free}.

For any graph $G$, consider the set of minimum semi-total dominating sets of vertices, $\{U_1,\dots,U_k\}$, and for $1\le i \le k$ let $X_i$ and $Y_i$ be partitions of $U_i$ into allied and free sets, respectively. We call $U_i$ so that $|X_i|$ is of maximum size for $1\le i \le k$ a \emph{maximum allied semi-total dominating set} of $G$, the partition $\{X_i,Y_i\}$ a \emph{maximum allied partition} of $G$, the set $X_i$ a \emph{maximum allied set} of $G$, and the set $Y_i$ a \emph{minimum free set} of $G$. 

For any maximum allied partition of $G$, $\{X,Y\}$, let $x(G)=|X|$ and $y(G)=|Y|$.

\begin{thm}\label{VizTypSemitot}
For any isolate-free graphs $G$ and $H$, \[\gamma_{t2}(G\cart H)\ge \frac{1}{3}\gamma_{t2}(G)\gamma_{t2}(H)\]
\end{thm}

\proof
\;Let $D$ be a minimum semi-total dominating set of $G\cart H$. Let $k=\gamma_{t2}(G)$ and $U=\{u_1,\dots, u_k\}$ be a maximum allied semi-total dominating set of $G$ with maximum allied partition $\{X,Y\}$. Suppose $X=\{u_1,\dots, u_{\ell}\}$ and $Y=\{u_{\ell+1},\dots, u_{\ell+m}\}$.

Form a partition $\{\pi_1,\dots, \pi_{\ell},\pi_{\ell+1},\dots, \pi_{\ell+m}\}$ of $V(G)$ where $\pi_i\subseteq N(u_i)$ and $x\in \pi_i$ implies $x$ is adjacent to $u_i$ for $1\le i \le \ell$, $\pi_j\subseteq N[u_j]$ and $x\in \pi_j$ implies $x$ is adjacent to $u_j$ for $\ell+1\le j\le \ell+m$. Furthermore, we define this partition to have the property that if $u_i\in X$ and $u_j\in Y$ so that $d(u_i,u_j)=2$, then $N(u_i)\cap N(u_j)\cap \pi_j=\emptyset$. That is, for any vertex $u_j$ of $Y$ which is of distance $2$ to some vertex of $X$, there exists a vertex $u$ which is adjacent to $u_j$ and to a vertex of $X$, and $u$ belongs to $\pi_i$ for some $i\in[\ell]$.

Let $D_i=(\pi_i\times V(H))\cap D$. Let $P_i=\{v:(u,v)\in D_i \mbox{ for some } u\in \pi_i\}$, which are the projections of $D_i$ onto $H$. We call vertices of $V(H)$ \emph{missing}, if they are not dominated from $P_i$ and write $M_i=V(H)-N_H[P_i]$. Vertices of $P_i$ which are of distance at most $2$ to some other vertex of $P_i$ or $M_i$ we call \emph{covered} and write $Q_i=\{v\in P_i : \exists w\in P_i\cup M_i \mbox{ such that } 0<d(v,w)\le 2\}$. Vertices of $P_i$ of distance at least $3$ to other vertices of $P_i$ or $M_i$ we call \emph{uncovered} and write $R_i=\{v\in P_i: \forall w\in (P_i\cup M_i)\backslash \{v\}, d(v,w)\ge 3\}$.

For $v\in V(H)$, let
\[D^v=D\cap (V(G)\times \{v\}) = \{(u,v)\in D:u\in V(G)\}\]
and $C$ be a subset of $\{1,\dots, k\}\times V(H)$ given by
\[C=\{(i,v):\pi_i\times \{v\} \subseteq N_{G \cart H}(D^v) \mbox{ or } v\in R_i\}.\]

Let $N=|C|$. We will bound $N$ from above by considering the following.
\[\cL_i=\{(i,v)\in C:v\in V(H)\},\]
\[\R^v=\{(i,v)\in C: 1\le i \le k\}.\]

These definitions are well-known as they appeared in the seminal work \cite{CS}, nonetheless, we would like to remind the reader of their interpretation. The set $C$ is a double indexing set, which indicates where you have cells that are either horizontally dominated or dominated by vertices of $R_i$. A cell is just a copy of $\pi_i$ for some $i$, at some height $v\in V(H)$. We represent $G$ along the horizontal axis of the Cartesian product and $H$ along the vertical. Thus, horizontally dominated cells are precisely, $\pi_i \times \{v\}$ which is contained in $N_{G\square H}(D^v)$. Now, $L_i$ are elements of $C$ with a fixed $i$ and $R^v$ are elements of $C$ along a fixed $v$.

Since counting vertices vertically and horizontally produces the same amount, we have
\[N=\sum_{i=1}^k|\cL_i|=\sum_{v\in V(H)}|\R^v|.\]

Notice that if $v\in M_i$, then the vertices in $\pi_i\times \{v\}$ which are not in $D^v$ must be adjacent to the vertices in $D^v$ since $D$ is a semi-total dominating set of $G\cart H$. Furthermore, the vertices of $R_i$ are counted in $\cL_i$. This means that $|\cL_i|\ge |M_i|+|R_i|$. Hence we obtain the following lower bound for $N$,
\begin{align}
N\ge \sum_{i=1}^k(|M_i|+|R_i|)\label{eq1}
\end{align}

To find an upper bound on the above quantity, we bound the size of $\R^v$.

\begin{claim}\label{upperboundN}
For any $v\in V(H)$, $|\R^v|\le 2|D^v|$.
\end{claim}
\proof
\;Suppose $|\R^v| > 2|D^v|$ for some $v\in V(H)$. For $(i,v)\in \R^v$, by definition, $\pi_i\times \{v\} \subseteq N_{G\cart H}(D^v)$ or $v\in R_i$. 

In what follows, we construct a semi-total dominating set $T$ of $G$.

 In the {\bf first case}, if $\pi_i\times \{v\} \subseteq N_{G\cart H}(D^v)$, we note that if some vertex $x\in \pi_i$, then $x$ is adjacent to vertices in $B^v$ where $B^v$ is the projection of $D^v$ onto $G$. 

\underline{Subcase 1.} {\it Suppose $u\in B^v$.} If $u\in \pi_i$ such that $(i,v)\notin \R^v$, $u\neq u_i$ and $1\le i \le \ell+m$, then $u\in N(u_i)$. If $u\in \pi_i$ such that $(i,v)\in \R^v$, then there exists $(u',v)\in B^v$ such that $u\in N(u')$. If $u\in \pi_i$ such that $(i,v)\notin \R^v$, $u=u_i$ for some $\ell+1\le i \le \ell+m$, then notice that we can find a vertex $x_i$ which is a neighbor of $u_i$ in $\pi_i$. Note that $x_i$ need not be a member of $B^v$, but simply a neighbor of $u_i$. Select one such vertex $x_i$ for every such $u$, and let $A$ be the set of these vertices $x_i$. Thus, $B^v\subseteq T$, $A:=\{u_i\,:\, (i,v)\notin \R^v,u_i\notin B^v,1\le i\le \ell+m\}\subseteq T$, and $\{x_i\,:\, (i,v)\notin\R^v, x_i\sim u \text{ for some } u\in U\cap B^v\}\subseteq T$.

\underline{Subcase 2.} {\it Suppose $u\in \{u_i:(i,v)\notin \R^v, 1\le i \le \ell\}$.} If $u\in \pi_j$ such that $(j,v)\notin \R^v$, then $u\in N(u_j)$. If $u\in \pi_j$ such that $(j,v)\in \R^v$, then there exists $(u',v)\in B^v$ such that $u\in N(u')$. Thus, in this subcase, $u$ is adjacent either to a vertex of $B^v$ or a vertex $u_j$. There are no new vertices that need to be added to $T$.

\underline{Subcase 3.} {\it Suppose $u\in \{u_i:(i,v)\notin \R^v, \ell+1\le i \le \ell+m\}$.} Suppose $u$ is of distance $2$ to some vertex $u_j\in X$. By the definition of the partition, there exists some vertex $w$ adjacent to $u$ and $u_j$, so that $w\in \pi_{j'}$ for some $j'\in [\ell]$. If $(j',v)\in \R^v$, then there exists $u'\in B^v$ so that $u'\sim w \sim u$, which means that $u$ is of distance at least $2$ to some vertex of $B^v$. Since $T$ contains $B^v$, these vertices are already distance $2$ from another vertex in $T$.

We are left to consider the case when $u$ is of distance at least $3$ to any vertex of $X$. Since $U$ is a minimum semi-total dominating set of $G$, there exists some vertex $u_j\in Y$, so that $d(u,u_j)=2$. If $(j,v)\notin \R^v$, these vertices are already in $T$ so no action needs to be taken.
% \john{Didn't we already handle these vertices in the previous paragraph? I feel like this can be simplified based on what is dominating and what is creating the semi-total dominating set.}
% \elliot{We haven't dealt with these vertices yet. The ones in the previous paragraph are at most of distance 2 to a vertex in X while these are farther away. Also, we are selecting $u_j$ here, not $u$.}

If $(j,v)\in \R^v$, then there exists some vertex $u'\in B^v$ so that $u'\sim u_j$. We will select $u_j$ and place it in $T$ to make $T$ a semi-total dominating set of $G$. Notice that in this case, the number of such vertices $u_j$ is at most equal to $|D^v|$.
Let $S$ be the set of such vertices $u_j$, which are of distance $2$ to a vertex $u\in Y$ and at least of distance $3$ to any vertex of $X$.
Then $S$ will be a subset of the set $T$. This finishes Subcase 3.

%Note also that if $x\in B^v$, then $x\in \pi_j$ where $(j,v)\notin \R^v$. Thus, notice that since $u_i$ is dominated by some vertex $u\in B^v$, we must have $d(u,u_j)\le 2$. 

In the {\bf second case}, if $v\in R_i$, then since $D$ is a semi-total dominating set, there is some vertex $(u,v)\in (\pi_i\times \{v\})\cap D^v$ and $(w,v)\in (\pi_j\times\{v\})\cap D^v$, for some $j\in[k]$, so that $(u,v)$ is at most distance $2$ from $(w,v)$. 

Putting these cases together, we have the following disjoint union of sets
\begin{align}
T = B^v &\cup \{u_i:(i,v)\notin \R^v, 1\le i \le \ell\} \cup \{u_i:(i,v)\notin \R^v, u_i\notin B^v, \ell+1\le i \le \ell+m\}\nonumber\\
&\cup A \cup S\label{semi-totalform}
\end{align}
To show $T$ is a semi-total dominating set of $G$, it is enough to show that $T$ is a dominating set, since we showed in each subcase of the first case, and in the second case, that every vertex of $T$ is of distance at most $2$ to some other vertex of $T$. If a vertex $u$ is contained in $\pi_i$ for $(i,v)\in \R^v$, then $u$ is dominated by some vertex of $B^v$. If $(i,v)\notin \R^v$, then $u$ is dominated either by $\{u_i:(i,v)\notin \R^v, 1\le i \le \ell\}$, or $\{u_i:(i,v)\notin \R^v, u_i\notin B^v, \ell+1\le i \le \ell+m\}$, or $A$.
 
Furthermore, 
\[|T|=|B^v|+(\gamma_{t2}(G)-|\R^v|+|S|)=2|D^v|+(\gamma_{t2}(G)-|\R^v|)<\gamma_{t2}(G)\]
which is a contradiction.
% \john{Explain why this is a semi-total dominating set.}
\smallqed

Thus, by claim \ref{upperboundN},
\begin{align}
N=\sum_{v\in V(H)}|\R^v|\le\sum_{v\in V(H)}2|D^v|=2|D|\label{eq2}
\end{align}

We now show a semi-total dominating set of $H$ in terms of $M_i$.

\begin{claim}\label{covering}
For any $i\in[k]$, there exists a set $X_i$ of at most $|R_i|-1$ vertices of $V(H)$ so that $M_i\cup P_i \cup X_i$ is a semi-total dominating set of $H$.
\end{claim}
\proof
\;We first observe that $P_i\cup M_i$ is a dominating set of $H$ with the additional property that the vertices of $M_i$ dominate only themselves, not their neighbors. Thus, every vertex $x\in R_i$ must be either of distance $3$ to some vertex $y\in R_i$ or every vertex of distance $2$ from $x$ is a vertex of $M_i$. This holds since otherwise some vertex of distance $2$ from $x$ is not dominated by $P_i\cup M_i$. Furthermore, if $x\in R_i$ which is of distance $3$ to some vertex $y\in R_i$, then we may select one vertex $z$ on a path from $x$ to $y$ such that $z$ is of distance at most $2$ to both $x$ and $y$.

We now construct a semi-total dominating set of $H$, $T_i$, by including the vertices of $M_i$, the vertices of $P_i$ and vertices $X_i$ which are of distance at most $2$ to two vertices of $R_i$ which are themselves of distance three to each other. The minimum number of such vertices is at most $|R_i|-1$, which can be easily verified by induction on $|R_i|$, and the result follows.
\smallqed 

By Claim \ref{covering}, for each $i$, we can construct a semi-total dominating set of $H$, $T_i=M_i\cup R_i\cup Q_i\cup X_i$. This gives $|M_i|+|R_i|\ge \gamma_{t2}(H)-|X_i|-|Q_i|$. However, note that $X_i\cap Q_i=\emptyset$ and $|X_i|+|Q_i|\le|P_i|$. This implies that $|M_i|+|R_i|\ge \gamma_{t2}(H)-|P_i|$. Thus, we have

\begin{align}
\sum_{i=1}^k\Big(|M_i|+|R_i|\Big)&\ge \sum_{i=1}^k\Big(\gamma_{t2}(H)-|P_i|\Big)=\gamma_{t2}(G)\gamma_{t2}(H)-|D|\label{eq3}
\end{align}

Combining equations \eqref{eq1}, \eqref{eq2}, and \eqref{eq3} we obtain
\begin{align*}
&|D|\ge\frac{1}{3}\gamma_{t2}(G)\gamma_{t2}(H)
\end{align*}
\qed

\section{Conclusion}
In this paper we have proven two Vizing's like results on the semi-total domination number. Our main result shown in Theorem~\ref{VizTypSemitot} shows that for isolate-free graphs $G$ and $H$, it must be the case that $\gamma_{t2}(G\cart H)\ge \frac{1}{3}\gamma_{t2}(G)\gamma_{t2}(H)$. However, we do not believe this bound is sharp, and conjecture a stronger result. 

\begin{conj}
For any isolate-free graphs $G$ and $H$, 
\begin{equation*}
\gamma_{t2}(G\:\Box \:H)\ge \frac{1}{2}\gamma_{t2}(G)\gamma_{t2}(H).
\end{equation*}
\end{conj}

\end{document}